\documentclass[12pt]{amsart}
\usepackage{amssymb}
\usepackage{eucal}
\newcommand{\rr}{\mathbb{R}}
\newcommand{\p}{\partial}
\newcommand{\zz}{\mathbb{Z}}

\newtheorem{lemma}{Lemma}[section]

\newtheorem{theorem}{Theorem}[section]

\author{L. Escauriaza, C. E. Kenig, G. Ponce, and L. Vega}
\title[KdV Equations]
{On uniqueness properties of solutions of 
the $k$-generalized KdV equations}

\begin{document}

\maketitle

\numberwithin{equation}{section}

\section{Introduction}

 In this paper we study uniqueness properties of solutions of the $k$-generalized Korteweg-de Vries 
equations
\begin{equation}
\label{1.1}
\p_t u + \p_x^3 u + u^k\p_x u=0, \;\;\;\;(x,t)\in  \rr^2,\;\;\;\;\;k\in\mathbb Z^+.
\end{equation}
 
 Our goal is to obtain sufficient conditions on the behavior of the difference  $u_1-u_2$ 
 of two solutions $u_1,\,u_2$ of \eqref{1.1} at two different times 
 $t_0=0$ and $t_1=1$ which guarantee that $u_1\equiv u_2$.

   This kind of uniqueness results has been deduced under the assumption that
 the solutions coincide in a large sub-domain of $\rr$ at two different times.
 In \cite{Zh}  B. Zhang proved that if $u_1(x,t)$ is a solution of the KdV, i.e. $k=1$ in \eqref{1.1},
  such that 
 $$
u_1(x,t)=0,\;\;\;\;(x,t)\in(b,\infty)\times \{t_0,t_1\}\,(\text{or}\,(-\infty,b)\times \{t_0,t_1\}) ,\;\;\;b\in\rr,
$$
then $\;u_1\equiv 0$, (notice that $u_2\equiv 0$ is a solution of \eqref{1.1}). His proof was based on the inverse scattering method (IST). In \cite{KePoVe4} this result was extended to any pair of solutions $u_1,\,u_2$ to  the generalized KdV equation, which includes non-integrable models. In particular, 
if $u_1,\,u_2$ are solutions of \eqref{1.1} in an appropriate class with $u_1(x,t)=u_2(x,t)$, for $(x,t)\in(b,\infty)\times \{t_0,t_1\}$ (or $(-\infty,b)\times \{t_0,t_1\})$, 
then $\;u_1\equiv u_2$.
  
In \cite{Ro} L. Robbiano proved the following uniqueness result :
  Let $u$ be a solution of the equation
 \begin{equation}
\label{1.2}
\p_t u + \p_x^3 u + a_2(x,t)\p^2_x u+a_1(x,t)\p_x u+a_0(x,t) u=0,   
\end{equation}
with coefficients $a_j,\,j=0,1,2$ in suitable function spaces. If $u(x,0)=0$ for $x\in(b,\infty)$ for some $b>0$ and there exist $c_1, \,c_2>0$ such that
$$
|\p_x^j u(x,t)|\leq c_1\,e^{-c_2x^{\alpha}}, \;\;\;\;\;\;\forall\;(x,t)\in (b,\infty)\times[0,1],\;\;\;j=0,1,2,
$$
for some $\alpha>9/4$, then $u\equiv 0$. This result applies to the difference
$u=u_1-u_2$ of two solutions $u_1,\,u_2$ of \eqref{1.1} with the  coefficients in \eqref{1.2} $a_0, a_1$ depending  on $u_1,\,u_2,\,\partial_x u_1,\,k$ and with $a_2\equiv 0$.

 In \cite{Ta}, using the IST,  S. Tarama showed that if the initial data $u(x,0)=u_0(x)$ has an appropriate exponential decay  for $x>0$, then the corresponding solution
of the KdV becomes analytic in the $x$-variable for all $t>0$. 

 It is interesting to notice that even  in the KdV case neither of the results  in \cite{Ro} and \cite{Ta} described above implies the other one. In \cite{Ro}
the decay assumption is needed in the whole time interval $[0,1]$, and the result in \cite{Ta} does not apply to the difference of two arbitrary solutions of the KdV.

Our main result concerning the equation \eqref{1.1} is the following.

  \begin{theorem}
 \label{Theorem 1.1}
 
 Let $u_1,\,u_2\in C([0,1]:H^2(\rr)\cap L^2(|x|^2dx))$ be strong solutions of \eqref{1.1}  in the domain
 $(x,t)\in\rr \times [0,1]$. If $k=1$ in \eqref{1.1} also assume that 
$u_1,\,u_2 \in C([0,1]:H^3(\rr))$. If 
 \begin{equation}
\label{1.3}
 u_1(\cdot,0)-u_2(\cdot,0),\;\;u_1(\cdot,1)-u_2(\cdot,1)\in H^1(e^{ax_{+}^{3/2}}dx),
 \end{equation}
 for any $a>0$,
 then $u_1\equiv u_2$. 
 
 \end{theorem}

 We shall say that  $f\in H^1(e^{ax_{+}^{3/2}}dx)$ if $f,\,\partial_{x} f\in L^2(e^{ax_{+}^{3/2}}dx)$,
 where $x_{+}=max\{x;\,0\}$, and $x_{-}=max\{-x;\,0\}$.

 \vskip.1in
 \underbar{Remarks}
 \vskip.1in
 a) The same result holds if in \eqref{1.3} instead of the space  $H^1(e^{ax_{+}^{3/2}}dx)$ 
 one considers  $H^1(e^{ax_{-}^{3/2}}dx)$.
 
 b) We recall that the solution of the associated linear initial value problem
 $$
 \p_t v + \p_x^3 v=0,\;\;\;\;\;\;\;\;\;v(x,0)=v_0(x),
 $$
 is given by the unitary group $\{U(t)\,:\,t\in \rr\}$ where 
 $$
 U(t)v_0(x)=\frac{1}{\sqrt[3]{3t}}\,Ai\left(\frac{\cdot}{\sqrt[3]{3t}}\right)\ast v_0(x), 
 $$
 and
 $$
 Ai(x)=\int_{\rr}\,e^{2\pi ix\xi+8\xi^3 i/3\pi^3}\,d\xi
 $$
 is the Airy function. This satisfies the estimate
 $$
 |Ai(x)|\leq c (1+x_{-})^{-1/4}\,e^{-c x_{+}^{3/2}}.
 $$
 
 Thus, the exponent $3/2$ in \eqref{1.3} can be seen as a reflection  of the asymptotic behavior of the Airy function. In fact, for the linear equation 
 $$
 \p_t v + \p_x^3  v=0,
 $$
 it shows that the decay rate in Theorem \ref{Theorem 1.1} is optimal.
 \vskip.1in
 c) In the particular case $u_2\equiv 0$ Theorem \ref{Theorem 1.1} tells us that 
 the only solution of the $k$-generalized KdV equation \eqref{1.1} which decays,  itself and its first derivative, 
 as $e^{-c x_{+}^{3/2}}$ at two different times is the zero solution.
 This is in contrast with the solutions of the equation
 $$
 \p_t u + \p_x^3 (u^2) + 2u \p_x u=0,
 $$mathbb
  found by Rosenau and Hyman  \cite{RoHy} called \lq\lq compactons". These are solitary waves of speed  $c$ with compact support 
 \begin{equation}
 \label{xxx}
 u_c(x,t)=
 \begin{cases}
 \begin{aligned}
&\tfrac{4c}{3}\,cos^2((x-ct)/4),\;\;\;\;\;\;\;\;\;\;\;\;\;|x-ct|\leq 2\pi,\\
&\;\;\;\;\;\;\;\;\;\;\;\;\;\;\;\;\;0,\;\;\;\;\;\;\;\;\;\;\;\;\;\;\;\;\;\;\;\;\;\;\;|x-ct|>2\pi.
\end{aligned}
\end{cases}
\end{equation}
  
d) In \cite{EsKePoVe} we proved the following result
 concerning the semi-linear Schr\"odinger equation
 \begin{equation}
\label{1.4}
i\p_t v + \Delta v + F(v,\overline v)=0,\;\;\;\;(x,t)\in\mathbb R^n\times \mathbb R.
\end{equation}

 \begin{theorem}
 \label{Theorem 1.2}
 
 Let $v_1,\,v_2\in C([0,1]:H^k(\mathbb R^n)),\,k\in \mathbb Z^+,\,k>n/2+1$ be strong solutions of the 
 equation \eqref{1.4} 
in the domain
 $(x,t)\in\mathbb R^n\times [0,1]$, with $F:\mathbb C^2\to \mathbb C$, $F\in C^k$  and $F(0)=\partial_uF(0)=\partial_{\bar u}F(0)=0$.
 If 
 \begin{equation}
\label{1.5}
 v_1(\cdot,0)-v_2(\cdot,0),\;\;v_1(\cdot,1)-v_2(\cdot,1)\in H^1(e^{a|x|^{2}}dx),
 \end{equation}
for any $a>0$,  then $v_1\equiv v_2$. 
 
 \end{theorem}
 
 The argument of the proof of Theorem 1.2 has  two main steps.
The first one is based on the exponential
 decay estimates obtained in \cite{KePoVe2}. These \lq\lq energy" estimates are  expressed in terms of the 
 $L^2(e^{\beta|x|}dx)$-norm and involve bounds independent of $\beta$. In \cite{EsKePoVe}
 they are used to deduce similar ones 
 with higher order powers in the exponent. The second step is to establish 
 lower bounds for the asymptotic behavior of the $L^2$-norm of the solution and its space gradient  in the annular  domain 
 $(x,t)\in \{R-1<|x|<R\}\times[0,1]$. This idea was motivated by the  work of Bourgain and Kenig \cite{BoKe} on a class of stationary Schr\"odinger operators ( i.e.  $-\Delta+V(x)$). Also in this second step we follow some arguments due to V. Izakov  \cite{Is}.
 
 For the equations  \eqref{1.1} considered here the first step in both \cite{KePoVe2} and \cite{EsKePoVe}, i.e.  weighted energy estimates, is not available. 
 We need to replace it by appropriate versions of Carleman estimates.  For example, for $H_{\beta}=(\p_t+e^{\beta x}\p_x^3e^{-\beta x})$ one has that
 \begin{equation}
\label{1.6}
\|e^{\beta x}\p_xe^{-\beta x} v\|_{L^{16}_xL^{16/5}_t}\leq c\|H_{\beta}v\|_{L^{16/15}_xL^{16/11}_t},
\end{equation}
 for functions $v\in C^{\infty}_0( \rr\times[0,1])$, see \cite{KePoVe4}. This kind of estimate resembles those  established  in \cite{KeRuSo} and  some extensions obtained in \cite{KePoVe4} related to the \lq\lq smoothing effect" found   in \cite{Ka},  \cite{KuFa} (homogeneous version)
 and in \cite{KePoVe1} (inhomogeneous version), see also \cite{IoKe}. However, we shall need their extension
 to functions $v\in C^{\infty}([0,1]:\mathbb{S}(\rr))$. In the case of \eqref{1.6} we shall prove that there exists $j\in\zz^+$ such that
 \begin{equation}
\label{1.7}
\|e^{\beta x}\p_xe^{-\beta x} v\|_{L^{16}_xL^{16/5}_t}\leq c \beta^j
(\| J^{1/2} v(\cdot,0)\|_{L^2}+
\| J^{1/2}v(\cdot,1)\|_{L^2})
+c\|H_{\beta}v\|_{L^{16/15}_xL^{16/11}_t}.
\end{equation}
 It will be crucial in our proof that  although in \eqref{1.7}  the constant in front of the norms involving the function $v$ evaluated  at 
 time $t=0$ and $t=1$ may grow as a power of $\beta$, the constant in front  of the norm of inhomogeneous term, i.e.  $H_{\beta}v$, is independent of $\beta>0$.

 e) Our argument here is direct and does not rely as that in \cite{KePoVe4}  on the unique continuation principle obtained by Saut
and Scheurer \cite{SaSc} :
if  a solution $v=v(x,t)$ of \eqref{1.2} 
in the domain $(x,t)\in (a,b)\times(t_1,t_2)$, with the coefficients $a_j,\,j=0,1,2$ in an appropriate class,
vanishes on an open set $\,\Omega\subseteq (a,b)\times(t_1,t_2)$, then
$v$ vanishes in the horizontal components of $\Omega$, i.e. the set
$$
\{(x,t)\in (a,b)\times(t_1,t_2)\,:\exists \,y\;\,s.t. \;(y,t)\in\Omega\}.
$$

 f) For the existence of solutions and well-posedness results for the IVP associated to the equation \eqref{1.1} we refer to \cite{KePoVe1} and references therein.
 We recall that the conditions $u_0\in H^2(\rr)\cap L^2(\rr:|x|^2dx)\equiv X_{2,2}$ and 
 $u_0\in  X_{2,2}\cap H^3(\rr)$  are locally preserved by the flow of solutions of  \eqref{1.1}, see \cite{Ka}, \cite{KePoVe1}. For our arguments  it suffices to have the decay 
 in only one side of the line, i.e. changing $|x|$ by $x_{+}$ in the weighted norms.
 This class is preserved for positive time $t>0$ by the flow of solutions, see \cite{KuFa}. 

 In particular, in the case $u_2\equiv 0$ we do not need any decay assumption on $u_1$
 since this will follow from the hypothesis \eqref{1.3}.
 
g)  Due to our interest  in results involving two different solutions $u_1,\,u_2$ of \eqref{1.1}
 we need to analyze  the equation satisfied by their difference $w=u_1-u_2$. This is  a linear equation
  of the form 
   \begin{equation}
 \label{1.8.b}
 \p_t w +\p_x^3 w + a_1(x,t)\p_x w+a_0(x,t) w=0,
 \end{equation}
 whose coefficients $a_1,\,a_0$ are polynomials of degree $k$ on $u_1,\,u_2$, and $\p_x u_1$. 
 Thus, the properties of $a_0,\,a_1$ depend
 on the class where the solutions $u_1,\,u_2$ are assumed and the value of $k$ considered.  
 
 In fact, we shall consider \eqref{1.2}, a more general equation than \eqref{1.8.b}.
 
  \begin{theorem}
 \label{Theorem 1.3}
 Assume that the coefficients in \eqref{1.2} satisfy that
 \begin{equation}
 \label{1.10}
\begin{aligned}
&a_0\in L^{4/3}_{xt}\cap L_x^{16/13}L^{16/9}_t\cap L_x^{8/7}L_t^{8/3},\\
&a_1\in L^{16/13}_xL^{16/9}_t\cap L_x^{8/7}L^{8/3}_t\cap L_x^{16/15}L_t^{16/3},\\
&a_2\in L^{8/7}_xL^{8/3}_t\cap L_x^{16/15}L^{16/3}_t\cap L_x^1L_t^{\infty}.\\
\end{aligned}
\end{equation}
Also, assume that  
\begin{equation}
\label{1.10.c}
a_0,\,a_1, \,a_2,\, \p_xa_2, \,\p_x^2 a_2\in L^{\infty}(\mathbb R^2),\,a_2,\, \p_ta_2\in L_t^{\infty}(\mathbb R:L^1_x(\mathbb R)).
\end{equation}

If $w\in C([0,1]:H^2(\rr)\cap L^2(|x|^2dx))$ is a  strong solution of \eqref{1.2}  in the domain
 $(x,t)\in\rr \times [0,1]$ with 
 \begin{equation}
\label{1.11}
 w(\cdot,0),\;w(\cdot,1)\in H^1(e^{ax_{+}^{3/2}}dx),
 \end{equation}
 for any $a>0$,
 then $w\equiv 0$. 
 
 \end{theorem}

 The remark (a) after the statement of Theorem 1.1 also applies here.
 
 As it was pointed out in the remark (b) the decay rate in \eqref{1.11}
 is optimal.
 
 We shall see that under the hypothesis of Theorem 1.1 the coefficients  $a_0,\;a_1$ of the equation 
 in \eqref{1.8.b}  belong to the class described in Theorem 1.3 in \eqref{1.10} and \eqref{1.10.c}. In fact,  it will be clear from our proof that the conditions in \eqref{1.10} in the $x$-variable are needed only in the positive semi-line,
 i.e. it suffices to have \eqref{1.10} with $a_j\in L^p_xL_t^q([0,\infty)\times [0,1])$ instead of
 $a_j\in L^p_xL_t^q=L^p_xL_t^q(\rr\times [0,1])$.
 
 It is here where the extra hypothesis $u_1,\,u_2\in C([0,1] : H^3)$ in Theorem 1.1 for the power $k=1$ in \eqref{1.1} is needed. 
    
 The rest of this paper is organized as follows. In section 2, we deduce 
 upper estimates in the time interval $[0,1]$ 
 for solutions of the inhomogeneous equation associated to \eqref{1.2} from the ones at times  $t_0=0$ and $t_1=1$
 and the inhomogeneous term. In 
 section 3, we shall obtain lower bounds
 for the $L^2$-norm of the solution and its first and second derivatives in  the annular domain mentioned above.
 Finally, in section 4 we combine the results in the previous  sections to prove Theorems \ref{Theorem 1.3} and \ref{Theorem 1.1}.
 
\section{Upper Estimates}

 We shall use the notations 
\begin{equation}
\label{2.1}
 Hf =(\p_t+\p_x^3)f,\;\;\;\;\;\;\;
 H_{\beta} f= ( \p_t + e^{\beta x} \p_x^3 e^{-\beta x} )f.
 \end{equation}

 Our first result in this section is the following lemma.
 
\begin{lemma}
\label{Lemma 2.1}

 There exists $k\in \zz^+$ such that if $ u\in C^{\infty}([0,1]:\mathbb{S}(\rr))$, then for any $\beta\geq 1$ 
\begin{equation}
\label{2.2}
\begin{aligned}
&\|e^{\beta x}u\|_{L^{8}_{xt}}+
\|e^{\beta x}\p_x u\|_{L^{16}_xL^{16/5}_t}+
\|e^{\beta x}\p_x^2 u\|_{L^{\infty}_xL^2_t}\\
&\leq c\beta^{2k} (\|J(e^{\beta x}u(\cdot,0))\|_{L^2}+
\|J(e^{\beta x}u(\cdot,1))\|_{L^2})\\
&+c(\|e^{\beta x}Hu\|_{L^{8/7}_{xt}}
+\|e^{\beta x}Hu\|_{L^{16/15}_xL^{16/11}_t}+
\|e^{\beta x}Hu\|_{L^1_xL^2_t}).
\end{aligned}
\end{equation}
where $\,Jg(x)=((1+|\xi|^2)^{1/2}\hat g)^{\lor}$ and  the norms  in the time variable 
(i.e. $\|\cdot \|_{L^p_t}$) are restricted  to the 
interval $[0,1]$.
\end{lemma}

In order to prove \eqref{2.2}, we set
\begin{equation}
\label{2.3}
v=e^{\beta x}u,
\end{equation}
and rewrite \eqref{2.2} as
\begin{equation}
\label{2.4}
\begin{aligned}
&
\|v\|_{L^{8}_{xt}}+
\|(e^{\beta x}\p_xe^{-\beta x})v\|_{L^{16}_xL^{16/5}_t}+
\|(e^{\beta x}\p^2_xe^{-\beta x})v\|_{L^{\infty}_xL^2_t}\\
&\leq c\beta^{2k} (\|J v(\cdot,0)\|_{L^2}+
\|J v(\cdot,1)\|_{L^2})\\
&+c(\|H_{\beta}v\|_{L^{8/7}_{xt}}
+\|H_{\beta}v\|_{L^{16/15}_xL^{16/11}_t}+
\|H_{\beta}v\|_{L^1_xL^2_t}).
\end{aligned}
\end{equation}

To obtain \eqref{2.4} we will prove the following string of inequalities
\begin{equation}
\label{2.6}
\|v\|_{L^8_{xt}}\leq c
(\| v(\cdot,0)\|_{L^2}+
\| v(\cdot,1)\|_{L^2})
+c\|H_{\beta}v\|_{L^{8/7}_{xt}},
\end{equation}
\begin{equation}
\label{2.7}
\|e^{\beta x}\p_xe^{-\beta x} v\|_{L^{16}_xL^{16/5}_t}\leq c \beta^k
(\| J^{1/2} v(\cdot,0)\|_{L^2}+
\| J^{1/2}v(\cdot,1)\|_{L^2})
+c\|H_{\beta}v\|_{L^{16/15}_xL^{16/11}_t},
\end{equation}
and
\begin{equation}
\label{2.8}
\|e^{\beta x}\p^2_xe^{-\beta x} v\|_{L^{\infty}_xL^2_t}\leq c \beta^{2k}
(\| J v(\cdot,0)\|_{L^2}+
\| J v(\cdot,1)\|_{L^2})
+c \|H_{\beta}v\|_{L^1_xL^2_t}.
\end{equation}

Clearly, \eqref{2.6}-\eqref{2.8} will imply \eqref{2.2}.

 First we shall prove first the following estimate which will be used later
\begin{equation}
\label{2.5}
\|v\|_{L^{\infty}_tL^2_x}\leq c
(\| v(\cdot,0)\|_{L^2}+
\| v(\cdot,1)\|_{L^2})
+c\|H_{\beta}v\|_{L^1_tL^2_x}.
\end{equation}

\it{Proof of \eqref{2.5}. }\rm  We have that
\begin{equation}
\label{2.9}
H_{\beta}=\p_t+e^{\beta x} \p_x^3 e^{-\beta x}=\p_t +(e^{\beta x} \p_x e^{-\beta x})^3,
\end{equation}
with
\begin{equation}
\label{2.10}
\begin{aligned}
&e^{\beta x} \p_x e^{-\beta x} =\p_x-\beta,\\
&(e^{\beta x} \p_x e^{-\beta x})^2 =(\p_x-\beta)^2=\p_x^2-2\beta \p_x+\beta^2,\\
&(e^{\beta x} \p_x e^{-\beta x})^3 =(\p_x-\beta)^3=\p_x^3-3\beta \p^2_x+3\beta^2\p_x-\beta^3\\
&=(\p_x^3+3\beta^2\p_x)-(3\beta\p_x^2+\beta^3)=\text{\it {skew-symmetric}}-\text{\it{symmetric part}}.
\end{aligned}
\end{equation}
The symbol of $H_{\beta}$ is
\begin{equation}
\label{2.11}
i\tau - i \xi^3+3 i \beta^2\xi -(\beta^3-3\beta\xi^2),
\end{equation}
whose real part $\beta^3-3\beta\xi^2$ vanishes at
\begin{equation}
\label{2.12}
\xi_{\pm}=\pm \beta/\sqrt{3},\;\;\;\;\;(\beta\geq 1).
\end{equation}
By an approximation argument it suffices to prove \eqref{2.5} for $v\in C^{\infty}([0,1]:\mathbb{S}(\rr))$
such that $\hat v(\xi,t)=0$ near  $\xi_{\pm}$ for all $t\in[0,1]$.

 Here, we shall denote by  $\,\hat f(\xi,t),\,\hat f(x,\tau),\,\hat f(\xi,\tau)$ the Fourier transform of $f(\cdot,\cdot)$ with respect to the dual variables $\xi,\;\tau,\,(\xi,\tau)$ respectively, i.e.  $\hat f(\cdot,\cdot)$ stands for  the Fourier transform of $f$ with respect to the dual variables where $\hat f$ is evaluated.

Assume now that $f\in C^{\infty}([0,1]:\mathbb{S}(\rr))$ with $f(x,t)= 0$ for $t$ near $0$ and $1$, so we can extend $f$ as $0$ outside the strip $\rr \times[0,1]$. Also assume  that $\hat f(\xi,t)=0$ for $\xi$ near 
$ \xi_{\pm}$ for all $t\in \rr$. Using our assumptions on $f$ we define
\begin{equation}
\label{2.13}
\widehat{Tf}(\xi,\tau)=\frac{\hat f(\xi,\tau)}{i\tau -i\xi^3+3i\beta^2\xi -(\beta^3-3\beta\xi^2)},
\end{equation}
 and claim that the estimate
 \begin{equation}
\label{2.14}
\|Tf\|_{L^{\infty}_tL^2_x}\leq c\|f\|_{L^1_tL^2_x}
\end{equation}
implies that in \eqref{2.5}.  To prove it we choose $\eta_{\epsilon}\in C^{\infty}(\rr)$, $\epsilon\in(0,1/4)$,  with
\begin{equation}
\label{2.15}
\eta_{\epsilon}(t)=1,\;\;\;t\in[2\epsilon,1-2\epsilon],\;\;\;\; \text{and}
\;\;\;\;\;supp\,\eta_{\epsilon}\subset[\epsilon,1-\epsilon],
\end{equation}
and define 
$$
v_{\epsilon}(x,t)=\eta_{\epsilon}(t)v(x,t),\;\;\;\;\;f_{\epsilon}(x,t)=H_{\beta}(v_{\epsilon})(x,t).
$$
Then, $v_{\epsilon}=Tf_{\epsilon}$ since both sides  have the same Fourier transform and both are
in $L^2_{xt}$  by our asssumptions on $ v$, which are inherited by $v_{\epsilon}$. Thus, \eqref{2.14}
gives
\begin{equation}
\label{2.16} 
\|v_{\epsilon}\|_{L^{\infty}_tL^2_x}\leq c\|H_{\beta}(v_{\epsilon})\|_{L^1_tL^2_x}
\leq c\|\eta'_{\epsilon}(t)v\|_{L^1_tL^2_x}+c\|\eta_{\epsilon}H_{\beta}(v)\|_{L^1_tL^2_x}.
\end{equation}
Letting $\epsilon\downarrow 0$ in \eqref{2.16} the left hand side converges to
$\|v\|_{L^{\infty}_{[0,1]}L^2_x}$, while  the right hand side has a limit  bounded by
$c(\|v_0\|_{L^2}+\|v_1\|_{L^2}) +c\|H_{\beta}v\|_{L^1_tL^2_x}$.

Hence, to obtain \eqref{2.5} we just need to prove \eqref{2.14}. In order to prove \eqref{2.14} it
suffices  to show that for $f(x,t)=f(x)\otimes \delta_{t_0}(t)$, 
with $\hat f(\xi)=0$ near $\pm\xi$, with $t_0\in(0,1)$ one has that

\begin{equation}
\label{2.17}
\|Tf\|_{L^{\infty}_tL^2_x}\leq c\|f\|_{L^2}, \;\;\;\;\text{with}\;\;c\;\;\text{independent of}\;\;t_0.
\end{equation}
First, we recall the formulas 
\begin{equation}
\label{2.18}
\left(\frac{1}{\tau+i b}\right)^{\vee}(t)=c
\begin{cases}
\begin{aligned}
&\chi_{(-\infty,0)}(t)\,e^{tb},\;\;\;\;b>0,\\
&\;\chi_{(0,\infty)}(t)\,e^{tb},\;\;\;\;b<0,
\end{aligned}
\end{cases}
\end{equation}
and consequently for $a,\,b\in \rr$
\begin{equation}
\label{2.19}
\left(\frac{e^{it_0\tau}}{\tau-a+i b}\right)^{\vee}(t)=c e^{ita}
\begin{cases}
\begin{aligned}
&\chi_{(-\infty,0)}(t-t_0)\,e^{(t-t_0)b},\;\;\;\;b>0,\\
&\;\chi_{(0,\infty)}(t-t_0)\,e^{(t-t_0)b},\;\;\;\;b<0.
\end{aligned}
\end{cases}
\end{equation}
Therefore,
\begin{equation}
\label{2.20}
\widehat{Tf}(\xi,\tau)=\frac{e^{it_0\tau} \hat f(\xi)}{i\{(\tau-\xi^3+3\beta^2\xi)+i(\beta^3-3\beta\xi^2)\}}
=-i\frac{e^{it_0\tau} \hat f(\xi)}{\tau-a(\xi)+ib(\xi)}.
\end{equation}
Combining \eqref{2.19}-\eqref{2.20} we see that the operator $T$ acting on these functions becomes the one variable operator $R$ 
\begin{equation}
\label{2.21}
\begin{aligned}
\widehat{Rf}(\xi)&=(\chi_{\{b(\xi)>0\}}(\xi)e^{ita(\xi)}e^{(t-t_0)b(\xi)}\chi_{(-\infty,0)}(t-t_0))\hat f(\xi)\\
&+(\chi_{\{b(\xi)<0\}}(\xi)e^{ita(\xi)}e^{(t-t_0)b(\xi)}\chi_{(0,\infty)}(t-t_0))\hat f(\xi),
\end{aligned}
\end{equation}
for which it needs to be established that
\begin{equation}
\label{2.22}
\|R f\|_{L^2_x}\leq c \|f\|_{L^2_x}, \;\;\;\text{with}\;\;c\;\;\text{independent of}\;\;\beta,\,t_0.
\end{equation}
But this  is obvious from the form of the multipliers in \eqref{2.21}. Therefore \eqref{2.5} is proved.
\vskip.1in
\it{Proof of \eqref{2.6}. }\rm  Again it suffices to show it for $v\in C^{\infty}([0,1]:\mathbb{S}(\rr))$
such that $\hat v(\xi,t)=0$ near $ \xi_{\pm}$.
Assume that $f\in C^{\infty}([0,1]:\mathbb{S}(\rr))$ with $f(x,t)= 0$ for $t$ near $0$ and $1$, so we can extend $f$ as $0$ outside the strip $\rr \times[0,1]$. Also assume that $\hat f(\xi,t)=0$ for $\xi$ near $ \xi_{\pm}$ for all $t\in \rr$. For the operator $T$ defined in \eqref{2.13} we shall show  for $f\in\mathbb{S}(\rr^2)$ with
$\hat f(\xi,\tau)=0$ near $\xi_{\pm}$ for all $t\in\rr$ 
{ \begin{equation}
\label{2.23}
\begin{aligned}
&(a)\;\;\;\|Tf\|_{L^8_{xt}}\leq c\|f\|_{L^{8/7}_{xt}},\\
&(b)\;\;\;\|Tf\|_{L^8_{xt}}\leq c\|f\|_{L^1_tL^2_x}.
\end{aligned}
\end{equation}
Assuming for the moment the inequalities in \eqref{2.23} we shall complete the proof of \eqref{2.6}. 

Consider 
\begin{equation}
\label{2.23a}
v_{\epsilon}(x,t)=\eta_{\epsilon}(t)v(x,t),\;\;\;\;\;H_{\beta}(v_{\epsilon})=\eta'_{\epsilon}(t)v+\eta_{\epsilon}H_{\beta}v=f_1(x,t)+f_2(x,t).
\end{equation}
Let 
\begin{equation}
\label{2.23b}
v_1(x,t)=Tf_1(x,t),\;\;\;\;\;\;v_2(x,t)=Tf_2(x,t),
\end{equation}
where both make sense by our assumption on $v$. Then,
\begin{equation}
\label{2.23c}
v_{\epsilon}(x,t)=v_1(x,t)+v_2(x,t),
\end{equation}
since both sides are in $L^2$ and have the same Fourier transform. Hence, from \eqref{2.23} it follows that
$$
\aligned
&\|v_{\epsilon}\|_{L^8_{xt}} \leq \|v_1\|_{L^8_{xt}}+\|v_2\|_{L^8_{xt}}\\
&\leq c \|f_1\|_{L^1_tL^2_x} + c \|f_2\|_{L^{8/7}_{xt}}
\leq c\|\eta'_{\epsilon}(t) v\|_{L^1_tL^2_x} +c\|\eta_{\epsilon}(t) H_{\beta}v\|_{L^{8/7}_{xt}},
\endaligned
$$
and letting $\epsilon \downarrow 0$ one gets \eqref{2.6}. So we need to establish \eqref{2.23}.
 \eqref{2.23}-(a) was proved in \cite{KePoVe4} (Lemma 2.3}). To obtain \eqref{2.23}-(b) we again restrict ourselves to consider 
 $f(x,t)=f(x)\otimes\delta_{t_0}(t)$, and reduce it to show that the operator $R$ defined in \eqref{2.21} 
 satisfies that
 $$
 \|Rf\|_{L^8_{xt}}\leq c\|f\|_{L^2},\;\;\;\text{with}\;\;c\;\;\text{independent of}\;\;\beta,\,t_0.
$$
But, this follows from the proof of Lemma 4.1 in \cite{KePoVe3}.

\vskip.1in
\it{Proof of \eqref{2.8}. }\rm  Again we make our usual assumptions on $\hat v(\xi, \tau)$. For $f\in\mathbb{S}(\rr^2)$ with $\hat f(\xi, t)=0$ near $\xi_{\pm}$ for all $t\in\rr$ we define using \eqref{2.13}
\begin{equation}
\label{2.24}
\widehat {T_2f} (\xi,\tau)=(i\xi-\beta)^2\,\widehat {Tf} (\xi,\tau)=\frac{(i\xi-\beta)^2\,\hat f(\xi,\tau)}{i\tau -i\xi^3+3i\beta^2\xi -(\beta^3-3\beta\xi^2)}.
\end{equation}
For the operator
\begin{equation}
\label{222}
\tilde  T_2f (x,t)=\chi_{[0,1]}(t) T_2f(x,t),
\end{equation}
we claim the following bounds
\begin{equation}
\label{2.25}
\begin{aligned}
(a) &\,\;\; \|\tilde  T_2 f \|_{L^{\infty}_xL^2_t}\leq c\|f\|_{L^1_xL^2_t},\\
(b) &\,\;\; \|\tilde T_2 f \|_{L^{\infty}_xL^2_t}\leq c \beta^2\| J f\|_{L^1_tL^2_x}.
\end{aligned}
\end{equation}
Assuming \eqref{2.25} (a)-(b) we shall prove \eqref{2.8}. With the notation in \eqref{2.23a}-\eqref{2.23c} 
from \eqref{2.25} it follows that
$$
\aligned
&\|(\p_x-\beta)^2v_{\epsilon}\|_{L^{\infty}_xL^2_t}
 \leq \|\chi_{[0,1]}(t)(\p_x-\beta)^2v_1\|_{L^{\infty}_xL^2_t}
 +\|\chi_{[0,1]}(t)(\p_x-\beta)^2v_2\|_{L^{\infty}_xL^2_t}\\
 &\leq  \|\tilde T_2f_1\|_{L^{\infty}_xL^2_t}+\|\tilde T_2f_2\|_{L^{\infty}_xL^2_t}
 \leq c\beta^2\| J f_1\|_{L^1_tL^2_x}+c\|f_2\|_{L^1_xL^2_t}\\
&\leq c\beta^2\|\eta'_{\epsilon}(t)Jv\|_{L^1_tL^2_x}+c\|\eta_{\epsilon}(t)H_{\beta}v\|_{L^1_xL^2_t},
\endaligned
$$
so letting $\epsilon\downarrow 0$ we obtain \eqref{2.8}.

The proof of
$$
\|T_2f\|_{L^\infty_xL^2_t}\leq c\|f\|_{L^1_xL^2_t},
$$
and therefore of  \eqref{2.25}-(a) follows with a minor modification 
from the argument given in \cite{KePoVe4}. Notice that the polynomial considered in the numerator of the 
fraction appearing  in (2.21) of  \cite{KePoVe4} is $\xi(i\xi -\beta)$
with $\beta=1$ while here we are considering $(i\xi -\beta)^2$. The  
proof works in exactly the same way as it can be easily checked. In fact, the use of this polynomial
instead of the one in  \cite{KePoVe4} is more convenient for the Littlewood-Paley interpolation 
argument which appears later on in (2.46)-(2.48) of  \cite{KePoVe4}. Notice that (2.48) 
in  \cite{KePoVe4} for $j=0$
is not true, but the proof just sketched fixed the error in  \cite{KePoVe4}. Another possible way to bypass this dfficulty is to use a Littlewood-Paley 
decomposition for $j\in \zz$ instead of $j=0,1,2...$.

We next prove \eqref{2.25}-(b). Let $\theta_r\in C^{\infty}_0(\rr)$ with $\theta_r(x)=1$, for $|x|\leq 2r$ and 
$supp\,\theta_r\subset\{|x|\leq 3r\}$, and consider
$$
\aligned
\widehat{T_2f}(\xi,\tau)&=
\frac{\theta_{\beta}(\xi)(i\xi-\beta)^2\,\hat f(\xi,\tau)}{i\tau -i\xi^3+3i\beta^2\xi -(\beta^3-3\beta\xi^2)}
+\frac{(1-\theta_{\beta}(\xi)(i\xi-\beta)^2\,\hat f(\xi,\tau)}{i\tau -i \xi^3+3i\beta^2\xi -(\beta^3-3\beta\xi^2)}\\
&=\widehat{T_{2,1}f}(\xi,\tau) + \widehat{T_{2,2}f}(\xi,\tau).
\endaligned
$$
Now, using Sobolev lemma one gets that
\begin{equation}
\label{2.26}
\|\tilde T_{2,1} f\|_{L^{\infty}_xL^2_t} \leq \|J \tilde T_{2,1}f\|_{L^2_xL^2_t}
=c\|J \tilde T_{2,1}f \|_{L^2_tL^2_x}\leq c\|JT_{2,1}f\|_{L^{\infty}_tL^2_x},
\end{equation}
where $\tilde T_{2,1}=\chi_{[0,1]}(t) T_{2,1}$. Now let
$$
\hat  g_1(\xi,\tau)=\theta_{\beta}(\xi)(1+|\xi|^2)^{1/2}(i\xi-\beta)^2\hat f(\xi,\tau).
$$
Then
$$
JT_{2,1}f(x,t)=Tg_1(x,t),
$$ 
so by \eqref{2.14} and \eqref{2.26} it follows that
\begin{equation}
\label{int1}
\|\tilde T_{2,1}f\|_{L^{\infty}_xL^2_t} \leq c\|g_1\|_{L^1_tL^2_x}
\leq c \beta^2\|J f\|_{L^1_tL^2_x}.
\end{equation}
To complete \eqref{2.25}-(b) it suffices to prove that
\begin{equation}
\label{int2}
\|T_{2,2}f\|_{L^{\infty}_xL^2_t}\leq c\|J f\|_{L^1_tL^2_x}.
\end{equation}
We again reduce ourselves to consider functions of the form 
$f(x,t)=f(x)\otimes\delta_{t_0}(t)$, so we just need to bound the operator
$$
\widehat{R_{2,2}f}(\xi,t)=
(1-\theta_{\beta}(\xi))(i\xi-\beta)^2\chi_{\{b(\xi)>0\}}(\xi)e^{ita(\xi)}e^{(t-t_0)b(\xi)}\chi_{(-\infty,0)}(t-t_0)\hat f(\xi)
$$
as
$$
\|R_{2,2}f\|_{L^{\infty}_xL^2_t}\leq c\|J f\|_{L^2_x},\;\text{with}\;\;c\;\;\text{independent of}\;\;\beta,\,t_0.
$$
We write
\begin{equation}
\label{2.27}
\begin{aligned}
&R_{2,2}f(x,t)=\\
&\int e^{ix\xi}(1-\theta_{\beta}(\xi))(i\xi-\beta)^2\chi_{\{b(\xi)>0\}}(\xi) 
e^{ita(\xi)}e^{(t-t_0)b(\xi)}\chi_{(-\infty,0)}(t-t_0)\hat f(\xi)d\xi,
\end{aligned}
\end{equation}
and recall that $a(\xi)=(\xi^3-3\beta^2\xi)$. Now we change variables
$$
\lambda=\xi^3-3\beta^2\xi,\;\;\;\;d\lambda=(3\xi^2-3\beta^2)d\xi=3(\xi^2-\beta^2)d\xi.
$$
From the definition of $\theta_{\beta}(\cdot)$ the domain of integration in \eqref{2.27} is equal to 
$\{|\xi|\geq 2\beta\}$, where $|\xi^2-\beta^2|\simeq |\xi|^2$, and the transformation is one to one.
Thus, we have $\xi=\xi(\lambda)$ and 
$$
\aligned
&R_{2,2}f(x,t)=\\
&\int e^{it\lambda} \;\frac{e^{ix\xi}(1-\theta_{\beta}(\xi))(i\xi-\beta)^2}{3(\xi^2-\beta^2)}\;\chi_{\{b(\xi)>0\}}(\xi)
e^{(t-t_0)b(\xi)}\chi_{(-\infty,0)}(t-t_0)\hat f(\xi)d\xi\\
&=\int e^{it\lambda}\hat g_2(\lambda)\Psi(\lambda,t)d \lambda,
\endaligned
$$
with
$$
\aligned
&\hat g_2(\lambda)=\frac{e^{ix\xi}(1-\theta_{\beta}(\xi))(i\xi-\beta)^2\hat f(\xi)}{3(\xi^2-\beta^2)},\\
&\Psi(\lambda,t)=\chi_{\{b(\xi)>0\}}(\xi) e^{(t-t_0)b(\xi)}\chi_{(-\infty,0)}(t-t_0).
\endaligned
$$
We observe that
$$
|\Psi(\lambda,t)|\leq c,\;\;\;\;\forall(\lambda,t)\in\rr^2\;\;\;\;\;\;\text{and}\;\;\;\;\;\int |\p_t\Psi(\lambda,t)|dt\leq c\;\;\;\;\;\forall \;\lambda\in\rr.
$$
Therefore, using the result in \cite{CoMe} (page 26) and taking adjoint one gets that
$$
\aligned
&\|\int e^{it\lambda}\hat g_2(\lambda)\Psi(\lambda,t)d \lambda\|_{L^2_t}\leq \|\hat g\|_{L^2}\\
&\leq c\left(\int\frac{|e^{ix\xi}(1-\theta_{\beta}(\xi))(i\xi-\beta)^2\hat f(\xi)|^2}{|3(\xi^2-\beta^2)|\,|3(\xi^2-\beta^2)|}
d\lambda\right)^{1/2}\\
&\leq c\left(\int\frac{|(1-\theta_{\beta}(\xi)|^2|\xi^2+\beta^2|^2|\hat f(\xi)|^2}{|3\xi^2-\beta^2|}d\xi\right)^{1/2}\\
&\leq c  \|Jf\|_{L^2},
\endaligned
$$
which finishes the proof of \eqref{2.25}-(b)


\vskip.1in
\it{Proof of \eqref{2.7}. }\rm  We make the usual assumptions on $v$ and $\hat v$. For $f\in\mathbb{S}(\rr^2)$ with $\hat f(\xi, t)=0$ near $\xi_{\pm}$ for all $t\in\rr$ we define using \eqref{2.13}
\begin{equation}
\label{2.28}
\widehat {T_1f} (\xi,\tau)=(i\xi-\beta)\,\widehat {Tf} (\xi,\tau)=\frac{(i\xi-\beta)\,\hat f(\xi,\tau)}{i\tau -i\xi^3+3i\beta^2\xi -(\beta^3-3\beta\xi^2)}.
\end{equation}
For the operator
\begin{equation}
\label{tilde}
\tilde  T_1f (x,t)=\chi_{[0,1]}(t) T_1f(x,t),
\end{equation}
we claim the following bounds
\begin{equation}
\label{2.29}
\begin{aligned}
(a) &\,\;\; \|\tilde  T_1 f \|_{L^{16}_xL^{16/5}_t}\leq c\|f\|_{L^{16/15}_xL^{16/11}_t},\\
(b) &\,\;\; \|\tilde T_1 f \|_{L^{16}_xL^{16/5}_t}\leq c \beta\| J^{1/2} f\|_{L^1_tL^2_x}.
\end{aligned}
\end{equation}
As above it  is easy to see that \eqref{2.7} follows from \eqref{2.29}.
Next, we recall that in \cite{KePoVe4} (see also the second paragraph after  \eqref{2.25} ) it was proved that
$$
\| T_1 f \|_{L^{16}_xL^{16/5}_t}\leq c\|f\|_{L^{16/15}_xL^{16/11}_t},
$$
which implies \eqref{2.29}-(a). To obtain \eqref{2.29}-(b) we write
$$
\aligned
\widehat{T_1f}(\xi,\tau)&=
\frac{\theta_{\beta}(\xi)(i\xi-\beta)\,\hat f(\xi,\tau)}{i\tau -i\xi^3+3i\beta^2\xi -(\beta^3-3\beta\xi^2)}
+\frac{(1-\theta_{\beta}(\xi)(i\xi-\beta)\,\hat f(\xi,\tau)}{i\tau -i \xi^3+3i\beta^2\xi -(\beta^3-3\beta\xi^2)}\\
&=\widehat{T_{1,1}f}(\xi,\tau) + \widehat{T_{1,2}f}(\xi,\tau),
\endaligned
$$
and consider first $\tilde T_{1,1}=\chi_{[0,1]}(t) T_{1,1}$. From \eqref{int1} 
$$
\|\tilde T_{2,1}f\|_{L^{\infty}_xL^2_t} 
\leq c \beta^2\|J f\|_{L^1_tL^2_x},
$$
and from \eqref{2.23}-(b) it follows that
$$
\|\tilde T_{0,1}f\|_{L^8_{xt}} 
\leq c  \|J f\|_{L^1_tL^2_x}.
$$
Hence using the interpolation argument based on the Littlewood-Paley decomposition as in
(2.46)-(2.48) of \cite{KePoVe4} one gets
$$
\|\tilde T_{1,1}f\|_{L^{16}_xL^{16/5}_t} 
\leq c \beta\|J^{1/2} f\|_{L^1_tL^2_x}.
$$
Finally, we interpolate between
$$
\|\tilde T_{0,2}f\|_{L^8_{xt}} 
\leq c \|J f\|_{L^1_tL^2_x},
$$
which follows from \eqref{2.23}-(a), with \eqref{int2} to get that
$$
\|\tilde T_{1,2}f\|_{L^{16}_xL^{16/5}_t} 
\leq c  \|J^{1/2} f\|_{L^1_tL^2_x},
$$
which yields \eqref{2.29}-(b).

This finished the proof of Lemma 2.1.

\vskip.1in 

Our next goal is to extend the estimates \eqref{2.2} in Lemma 2.1 to solutions of the linear equation with variable coefficients
\begin{equation}
\label{2.30}
\p_t u+\p_x^3 u +a_2(x,t)\p_x^2 u + a_1(x,t) \p_x u +a_0(x,t) u =0.
\end{equation} 
We introduce the notation
\begin{equation}
\label{2.31}
H_a=\p_t +\p_x^3 +a_2(x,t)\p_x^2  + a_1(x,t) \p_x  +a_0(x,t),
\end{equation} 
and try to find conditions which guarantee that multiplication by 
$a_0(x,t)$ maps
\begin{equation}
\label{2.32}
L^8_{xt}\to L^{8/7}_{xt},\;\;\;\;\;\;\;L^8_{xt}\to L^{16/15}_xL^{16/11}_t,\;\;\;\;\;\;\;L^8_{xt}\to L^1_xL^2_t,
\end{equation}
multiplication by 
$a_1(x,t)$ maps
\begin{equation}
\label{2.33}
L^{16}_xL^{16/5}_t\to L^{8/7}_{xt},\;\;\;\;\;\;L^{16}_xL^{16/5}_t\to L^{16/15}_xL^{16/11}_t,\;\;\;\;\;L^{16}_xL^{16/5}_t\to L^1_xL^2_t,
\end{equation}
and multiplication by 
$a_2(x,t)$ maps
\begin{equation}
\label{2.34}
L^{\infty}_xL^{2}_t\to L^{8/7}_{xt},\;\;\;\;\;\;L^{16}_xL^{16/5}_t\to L^{16/15}_xL^{16/11}_t,\;\;\;\;\;\;\;L^{16}_xL^{16/5}_t\to L^1_xL^2_t.
\end{equation}

So it suffices to have that
\begin{equation}
\begin{aligned}
\label{3-spaces}
&a_0\in L^{4/3}_{xt}\cap L_x^{16/13}L^{16/9}_t\cap L_x^{8/7}L_t^{8/3},\\
&a_1\in L^{16/13}_xL^{16/9}_t\cap L_x^{8/7}L^{8/3}_t\cap L_x^{16/15}L_t^{16/3},\\
&a_2\in L^{8/7}_xL^{8/3}_t\cap L_x^{16/15}L^{16/3}_t\cap L_x^1L_t^{\infty}.\\
\end{aligned}
\end{equation}
Thus, if $a_0,\,a_1,\,a_2$ are in these spaces with small norm, then the  inequality \eqref{2.2} will hold with $H_a$ in \eqref{2.31} instead of $H$, and one has the following result.

\begin{lemma}
\label{Lemma 2.2}

Assume that the coefficients in \eqref{2.30} $a_0,\,a_1,\,a_2$ satisfy \eqref{3-spaces} with small enough norms. There exists $k\in \zz^+$ such that if $ u\in C^{\infty}([0,1]:\mathbb{S}(\rr))$, then for any $\beta\geq 1$ 
\begin{equation}
\label{2.2.a}
\begin{aligned}
&\|e^{\beta x}u\|_{L^{8}_{xt}}+
\|e^{\beta x}\p_x u\|_{L^{16}_xL^{16/5}_t}+
\|e^{\beta x}\p_x^2 u\|_{L^{\infty}_xL^2_t}\\
&\leq c\beta^{2k} (\|J(e^{\beta x}u(\cdot,0))\|_{L^2}+
\|J(e^{\beta x}u(\cdot,1))\|_{L^2})\\
&+c(\|e^{\beta x}H_au\|_{L^{8/7}_{xt}}
+\|e^{\beta x}H_au\|_{L^{16/15}_xL^{16/11}_t}+
\|e^{\beta x}H_au\|_{L^1_xL^2_t}).
\end{aligned}
\end{equation}
\end{lemma}

\begin{proof}
\label{proof of Lemma 2.2}

First we introduce the notation
\begin{equation}
\begin{aligned}
\label{norms}
&|||h|||_1\equiv \|e^{\beta x} h\|_{L^{8}_{xt}}+\|e^{\beta x} \p_x h\|_{L^{16}_xL^{16/5}_t}+\|e^{\beta x} \p_x^2h\|_{L^{\infty}_xL^2_t},\\
&|||h|||_2\equiv \|h\|_{L^{8/7}_{xt}}+\|h\|_{L^{16/15}_xL^{16/11}_t}+\|h\|_{L^1_xL^2_t},
\end{aligned}
\end{equation}
From Lemma 2.1 and our assumptions it follows that
\begin{equation}
\label{2.34.a}
\begin{aligned}
||| u|||_1&\leq c\beta^{2k} (\|J(e^{\beta x}u(\cdot,0))\|_{L^2}+
\|J(e^{\beta x}u(\cdot,1))\|_{L^2})\\
&+||| e^{\beta x} H u|||_2\\
& \leq c\beta^{2k} (\|J(e^{\beta x}u(\cdot,0))\|_{L^2}+
\|J(e^{\beta x}u(\cdot,1))\|_{L^2})\\
&+||| e^{\beta x} H_a u|||_2+||| e^{\beta x} (a_2(x,t) \p_x^2  + a_1(x,t) \p_x  +a_0(x,t) )u|||_2\\
&\leq c\beta^{2k} (\|J(e^{\beta x}u(\cdot,0))\|_{L^2}+
\|J(e^{\beta x}u(\cdot,1))\|_{L^2})\\
&+||| e^{\beta x} H_a u|||_2+\frac{1}{2} \,||| u|||_1. 
\end{aligned}
\end{equation}
Hence,
\begin{equation}
\label{2.35.a}
|||e^{\beta x}u|||_1 
\leq c\beta^{2k} (\|J(e^{\beta x}u(0))\|_{L^2}+\|J(e^{\beta x}u(1))\|_{L^2})+c |||e^{\beta x} H_a u|||_2,
\end{equation}
which yields the desired result.
\end{proof}

 We now start with $u$ solving 
\begin{equation}
\label{equ}
 \p_tu +\p_x^3 u + a_2(x,t)\p_x^2 u + a_1(x,t) \p_x u +a_0(x,t) u =0,\;\;\;\;\;\;(x,t)\in\rr\times [0,1],
\end{equation}
with $u_0=u(\cdot,0),\,u_1=u(\cdot,1)\in H^1(e^{ax_{+}^{\alpha}})$ for some $a>0$, $\alpha>1$, and $a_0,\,a_1,\,a_2$ just in the spaces
in \eqref{3-spaces}. 

Choose $R$ so large that in the $x$-interval $(R,\infty)$  
the coefficients  $a_0,\,a_1,\,a_2$ in the corresponding spaces   \eqref{3-spaces} have small
norms.  Let $\mu\in C^{\infty}(\rr)$ with  $\mu(x)=0$ if $x<1$, and $\mu(x)=1$ if $x>2$. For 
$\mu_R(x)=\mu(x/R)$ we have that 
$$
u_R(x,t)=\mu_R(x) u(x,t),
$$
satisfies the equation
$$
 \p_tu_R +\p_x^3 u_R + a_2(x,t)\p_x^2 u_R + a_1(x,t) \p_x u_R +a_0(x,t) u_R = e_R(x,t),
 $$
where
$$
\aligned
&e_R(x,t)= \mu_R^{(3)}\tfrac{1}{R^3} u + 
3 \mu_R^{(2)}\tfrac{1}{R^2}\p_xu
+3 \mu_R^{(1)}\tfrac{1}{R} \p_x^2u \\
&+ a_2(x,t)(2 \mu_R^{(1)}\tfrac{1}{R} \p_xu+
\mu_R^{(2)}\tfrac{1}{R^2} u) +
a_1(x,t) \mu_R^{(1)}\tfrac{1}{R} u.
\endaligned
$$

Notice that $\,supp \;e_R\subset\{x\,:\,R<x<2R\}$. We will take 
$$
\beta=\frac{a}{2}\,R^{(\alpha-1)}.
$$
 Now we apply our inequality \eqref{2.2.a} to $u_R$, with 
 $$
 H_{a\tilde \mu_R}=\p_t +\p_x^3 +a_2(x,t)\tilde \mu_R(x) \p_x^2  + a_1(x,t) \tilde \mu_R(x) \p_x  +a_0(x,t) \tilde \mu_R(x).
$$
where $\tilde \mu_R(x) \mu_R(x)=\mu_R(x)$, and so that $a_j(x,t)\tilde \mu_R(x)$ with $j=0,1,2$ have small norm in the corresponding  spaces in  \eqref{3-spaces}  for $R>R_0$. From Lemma 2.2 it follows that  for $R$ large 
\begin{equation}
\label{2.34.c}
||| u_R|||_1\leq c\beta^{2k} (\|J(e^{\beta x}u_R(\cdot,0))\|_{L^2}+
\|J(e^{\beta x}u_R(\cdot,1))\|_{L^2})
 + |||e^{\beta x} e_R|||_2. 
\end{equation}

To bound the first term in the right hand side of \eqref{2.34.c} we use that
\begin{equation}
\label{2.36.c}
\begin{aligned}
&(1+\beta^{2k})\|J(e^{\beta x}u_R(0))\|_{L^2}
\leq c(1+\beta^{2k+1})(\| e^{\beta x} u_R(0)\|_{L^2}+\| e^{\beta x}\partial_xu_R(0)\|_{L^2})\\
&\leq c(1+\beta^{2k+1})( \| e^{\beta x}u(0)\|_{L^2\{x>R\}} + \|e^{\beta x}\p_xu(0)\|_{L^2\{x>R\}})\\
&\leq c(1+\beta^{2k+1})(\|e^{aR^{\alpha-1}x/2}u(0)\|_{L^2\{x>R\}} + \|e^{aR^{\alpha-1}x/2}\p_x u(0)\|_{L^2\{x>R\}}).
\end{aligned}
\end{equation}

 Since $k\in\zz^+$ is fixed and $\beta=a\,R^{(\alpha-1)}/2$ for $R$ sufficiently large,
 depending on $\alpha$ and $a$, one has 
\begin{equation}
\label{constant}
a^{2k+1} R^{(2k+1)(\alpha-1)} e^{aR^{\alpha-1}x/2} \leq c_{a,\alpha} e^{ax^{\alpha}}, \;\;\;\;\text{for}\;\;x>R>0.
\end{equation}
Then the right hand side of \eqref{2.36.c} is bounded by $c_{a,\alpha}$. 

A similar argument shows that 
$$
(1+\beta^{2k})\|J(e^{\beta x}u_R(1))\|_{L^2}\leq c_{a,\alpha}.
$$

Thus, it remains to bound $|||e^{\beta x} e_R|||_2$. Since $\,supp \;e_R\subset\{x\,:\,R<x<2R\}$ combining H\"older inequality and Minkowski's integral inequality it follows that
$$
|||e^{\beta x}e_R|||_2\leq c e^{aR^{\alpha-1} R}\,\| (|u| +|\p_x u| +|\p_x^2 u| )\,\chi_{\{x\,:\,R<x<2R\}}\|_{L^{\infty}_tL^2_x}\leq c' \,e^{aR^{\alpha-1} R}.
$$

Inserting these estimates in \eqref{2.34.c} we obtain that
$$
\aligned
&\|e^{\beta x} u\|_{L^{8}_{\{x>4R\}}L^8_t}+\|e^{\beta x} \p_xu\|_{L^{16}_{\{x>4R\}} L^{16/5}_t}
+\|e^{\beta x} \p_x^2 u\|_{L^{\infty}_{\{x>4R\}}L^2_t}\\
&\leq c_{a,\alpha} + c'  e^{aR^{\alpha-1} R}=c_{a,\alpha} + c'  e^{aR^{\alpha}}.
\endaligned
$$

If $\,x>4R$, then $\,e^{aR^{\alpha-1} x/2}\, e^{-aR^{\alpha}} \geq e^{aR^{\alpha}} $, so we get
\begin{equation}
\label{2.36.b}
\begin{aligned}
&e^{aR^{\alpha}} (\|u\|_{L^{8}_{\{x>4R\}}L^8_t}+\|\p_x u\|_{L^{16}_{\{x>4R\}} L^{16/5}_t}
+\|\p_x^2 u\|_{L^{\infty}_{\{x>4R\}}L^2_t})\\
& \leq 
 c_{\alpha,a}.
\end{aligned}
\end{equation}
 Therefore, using H\"older inequality in \eqref{2.36.b} it follows that for $R$ sufficiently large
 $$
 \aligned
 &\|u\|_{L^2(\{4R<x<4R+1\}\times (0,1))}+
  \|\p_xu\|_{L^2(\{4R<x<4R+1\}\times (0,1))}\\
  &+ \|\p_x^2 u\|_{L^2(\{4R<x<4R+1\}\times (0,1))}\leq c_{a,\alpha} e^{-aR^{\alpha}}.
  \endaligned
 $$
Now changing $4R$ by $R'$ we get that for any $R'>0$ sufficiently large
$$
\aligned
& \|u\|_{L^2(\{R'<x<R'+1\}\times (0,1))}+
  \|\p_xu\|_{L^2(\{R'<x<R'+1\}\times (0,1))} + \|\p_x^2 u\|_{L^2(\{R'<x<R'+1\}\times (0,1))}\\
  &\leq c_{a,\alpha} e^{-a(R')^{\alpha}/4^{\alpha}}.
 \endaligned
 $$
 So we have proved the following upper estimates for solutions of \eqref{2.30}.
 
 \begin{theorem}
\label{Theorem 2.3}

Assume that the coefficients in \eqref{equ} $a_0,\,a_1,\,a_2$ satisfy \eqref{3-spaces}. If $u=u(x,t)$ is a solution of \eqref{equ} 
with $u\in C([0,1]:H^2(\rr))$
satisfying that
$$
u(\cdot,0),\;\;\;u(\cdot,1)\in H^1(e^{ax_{+}^{\alpha}})
$$
 for some $\alpha>1$ and $a>0$, 
then there exist $c_0$ and $R_0>0$ sufficiently large such that for $R\geq R_0$ 
$$
\aligned
& \|u\|_{L^2(\{R<x<R+1\}\times (0,1))}+
  \|\p_xu\|_{L^2(\{R<x<R+1\}\times (0,1))}\\
  \\
  &+
   \|\p_x^2 u\|_{L^2(\{R<x<R+1\}\times (0,1))}\leq c_0 e^{-aR^{\alpha}/4^{\alpha}}.
 \endaligned
 $$

\end{theorem}
 
 \section{Lower Bounds}
 
 This section is concerned with lower bound estimates for the $L^2$-norm of a
 solution $u$ of the equation  \eqref{1.2} and its first order space derivative 
 $\p_xu$ in the box 
  $\{R-1<x<R\}\times[0,1]$. 
 
 \begin{lemma}
 \label{Lemma 3.1}
Assume that $\varphi:[0,1]\longrightarrow\rr$ is a smooth function. Then, there exist 
$c >0$ 
 and $M_1=M_1(\|\varphi'\|_{\infty};\|\varphi''\|_{\infty})>0$ such that the inequality
\begin{equation}
\label{3.1}
\begin{aligned}
&\frac{\alpha^{5/2}}{R^3}\, \|e^{\alpha (\frac xR+\varphi(t))^2} (\tfrac xR+\varphi(t))^2g\|_{L^2(dxdt)}
+
\frac{\alpha^{3/2}}{R^2}\, \|e^{\alpha (\frac xR+\varphi(t))^2} (\tfrac xR+\varphi(t))\p_xg\|_{L^2(dxdt)}\\
\\
&
+\frac{\alpha^{1/2}}{R}\, \|e^{\alpha (\frac xR+\varphi(t))^2} \p_x^2g\|_{L^2(dxdt)}
\,\leq \,c \,\| e^{\alpha (\frac xR+\varphi(t))^2}(\partial_t+\p_x^3)g\|_{L^2(dxdt)}
\end{aligned}
\end{equation}
holds, for $R\geq 1$, $\alpha$ such that $\alpha^2\geq M_1R^3$, and $g\in C_0^\infty(\rr^{2})$ supported  in 
$$
\{(x,t)\in\rr^2 : |\tfrac xR+\varphi(t)|\ge 1\}.
$$

\end{lemma}

\begin{proof}
\label{Proof of Lemma 3.1}

We define $f=e^{\alpha \theta(x,t)}g$,
for a general smooth function $\theta(x,t)$,
and consider the expression
\begin{equation}
\label{3.2}
e^{\alpha \theta(x,t)}(\p_t+\p_x^3)(e^{-\alpha \theta(x,t)}f(x,t))
=S_\alpha f + A_\alpha f\ ,
\end{equation}
where
\begin{align*}
&S_\alpha=-3\alpha\p_x(\p_x\theta(x,t)\p_x\cdot) \,+\, (-\alpha^3(\p_x \theta(x,t))^3-
\alpha\p_x^3\theta(x,t)-\alpha\p_t\theta(x,t))\\
\\
&A_\alpha=\p_t+\p_x^3+3\alpha^2(\p_x \theta(x,t))^2\p_x + 3\alpha^2 \p_x\theta(x,t) \p_x^2\theta(x,t).
\end{align*}
Thus,
\begin{equation}
\label{3.3}
A_\alpha^*=-A_\alpha,\;\;\;\;\;\;\;S_\alpha^*=S_\alpha,
\end{equation}
and one has
\begin{equation}
\label{3.4}
\begin{aligned}
&\| e^{\alpha \theta(x,t)}(\partial_t+\p_x^3)g\|^2_{L^2(dxdt)}= \|(A+S)f\|_{L^2(dxdt)}^2\\
&=\langle (A+S)f, (A+S)f\rangle\\
&=\|Af\|_2^2+\|Sf\|_2^2+\langle Af,Sf\rangle+\langle Sf,Af\rangle\\
& \geq \langle (SA-AS)f,f\rangle.
\end{aligned}
\end{equation}
 A computation shows that
\begin{equation}
\label{3.5}
\begin{aligned}
 &(SA-AS)f=[S;A]f=9\alpha\p_x^2(\p_x^2 \theta\,\p_x^2f)\\
 &+\p_x((6\alpha\p_x^4 \theta+6\alpha \p^2_{tx}\theta-18\alpha^3(\p_x \theta)^2\p_x^2\theta)\p_x f)\\
 &+(-3\alpha^3(\p_x^2\theta)^3-18\alpha^3\p_x\theta\,\p^2_x\theta\,
 \p_x^3\theta-3\alpha^3(\p_x\theta)^2\p_x^4\theta+\alpha\p_x^6\theta\\
 &+2\alpha\p^3_x\p_t\theta+\alpha\p_t^2\theta
 +6\alpha^3(\p_x\theta)^2\p_{tx}^2\theta+9\alpha^5(\p_x\theta)^4\p_x^2\theta)f.
 \end{aligned}
 \end{equation}

Now taking $\theta(x,t)=(x/R+\varphi(t))^2$ it follows from \eqref{3.5} and integrations by parts that
\begin{equation}
\label{3.6}
\begin{aligned}
&\langle (SA-AS)f,f\rangle =\frac{18\alpha}{R^2} \|\p_x^2 f\|_{L^2(dxdt)}^2\\
&-\frac{12\alpha}{R}\int\int \varphi'(t)(\p_xf)^2dxdt+\frac{144\alpha^3}{R^4}
\int\int (\tfrac xR+\varphi(t))^2(\p_xf)^2dxdt\\
&-\frac{24\alpha^3}{R^6}\int\int f^2dxdt+2\alpha\int\int (\varphi'(t)f)^2dxdt
+2\alpha\int\int (\tfrac xR+\varphi(t))\varphi''(t)f^2dxdt\\
&+\frac{48\alpha^3}{R^3}\int\int\varphi'(t)(\tfrac xR+\varphi(t))^2f^2dxdt
+\frac{288\alpha^5}{R^6}\int\int(\tfrac xR+\varphi(t))^4f^2dxdt\\
&=I_1+I_2+I_3+I_4+I_5+I_6+I_7+I_8.
\end{aligned}
\end{equation}

We first observe that 
$$
I_5+I_7+I_8=2\,\int\int\,(\alpha^{1/2}\varphi'(t)f+\tfrac{12\,\alpha^{5/2}}{R^3}(\tfrac xR+\varphi(t))^2\,f)^2dxdt.
$$
Therefore, since $|\tfrac xR+\varphi(t)|>1$ on the support of $f$, for 
\begin{equation}
\label{3.7}
\alpha^2\geq \|\varphi'\|_{\infty}R^3,
\end{equation}
it follows that 
$$
I_5+I_7+I_8\geq \frac{242\alpha^5}{R^6}\int\int(\tfrac xR+\varphi(t))^4f^2dxdt.
$$
Similarly, since $|\tfrac xR+\varphi(t)|>1$ on the support of $f$ for 
\begin{equation}
\label{3.8}
\alpha^2\geq (\|\varphi''\|^{1/2}_{\infty}+1) R^3,
\end{equation}
it follows that 
$$
\frac{2\alpha^5}{R^6}\int\int(\tfrac xR+\varphi(t))^4f^2dxdt\geq |I_6|,
$$
and that
$$
\frac{24\alpha^5}{R^6}\int\int(\tfrac xR+\varphi(t))^4f^2dxdt\geq |I_4|.
$$
Also from \eqref{3.7} one has that
$$
I_2+I_3\geq \frac{132\alpha^3}{R^4}
\int\int (\tfrac xR+\varphi(t))^2(\p_xf)^2dxdt.
$$

Hence, gathering the above information we conclude that
for 
\begin{equation}
\label{3.9}
\alpha^2\geq (\|\varphi'\|_{\infty}+\|\varphi''\|^{1/2}_{\infty}+1) R^3,
\end{equation}
one has that
\begin{equation}
\label{3.10}
\begin{aligned}
&\| e^{\alpha(\tfrac xR+\varphi(t))^2}(\partial_t+\p_x^3)g\|^2_{L^2(dxdt)}= \|(A+S)f\|_{L^2(dxdt)}^2\\
&\geq \langle (SA-AS)f,f\rangle
\geq \frac{18\alpha}{R^2} \int\int (\p_x^2 f)^2dxdt\\
&+\frac{132\alpha^3}{R^4}
\int\int (\tfrac xR+\varphi(t))^2(\p_xf)^2dxdt
+\frac{216\alpha^5}{R^6}\int\int(\tfrac xR+\varphi(t))^4f^2dxdt.
\end{aligned}
\end{equation}

Next, we rewrite \eqref{3.10} in terms of $g=e^{-\alpha(\tfrac xR+\varphi(t))^2 }f$.
In fact, it follows from \eqref{3.9} and \eqref{3.10} that there exits a universal constant $c_0>0$ such that
\begin{equation}
\label{3.11}
\begin{aligned}
\| e^{\alpha(\tfrac xR+\varphi(t))^2}&(\partial_t+\p_x^3)g\|_{L^2(dxdt)}
\geq \frac{c_0\alpha^{1/2}}{R}
 \left(\int\int e^{2\alpha(\tfrac xR+\varphi(t))^2 } (\p_x^2 g)^2dxdt\right)^{1/2}\\
&+\frac{c_0\alpha^{3/2}}{R^2}\left(
\int\int (\tfrac xR+\varphi(t))^2e^{2\alpha(\tfrac xR+\varphi(t))^2 }(\p_xg)^2dxdt\right)^{1/2}\\
&+\frac{c_0\alpha^{5/2}}{R^3}\left(\int\int(\tfrac xR+\varphi(t))^4e^{2\alpha(\tfrac xR+\varphi(t))^2 }g^2dxdt\right)^{1/2},
\end{aligned}
\end{equation}
which completes the proof of Lemma 3.1.
\end{proof}

Next, we shall extend the result of Lemma 3.1
to operators of the form
\begin{equation}
\label{3.12}
L=\p_t+\p_x^3 +a_0(x,t)+a_1(x,t)\p_x,
\end{equation}
with
$$
a_0,\;a_1\in L^{\infty}(\mathbb R^2).
$$

 \begin{lemma}
 \label{Lemma 3.2}
Assume that $\varphi:[0,1]\longrightarrow\rr$ is a smooth function. Then, there exist 
$c>0$, 
$ R_0=R_0( \|\varphi'\|_{\infty};\|\varphi''\|_{\infty};\|a_0\|_{\infty}; \|a_1\|_{\infty})>1$ and $M_1=M_1(\|\varphi'\|_{\infty};\|\varphi''\|_{\infty})>0$ such that the inequality
\begin{equation}
\label{3.13}
\begin{aligned}
&\frac{\alpha^{5/2}}{R^3}\, \|e^{\alpha (\frac xR+\varphi(t))^2} (\tfrac xR+\varphi(t))^2g\|_{L^2(dxdt)}
+\frac{\alpha^{3/2}}{R^2}\, \|e^{\alpha (\frac xR+\varphi(t))^2} \p_x g\|_{L^2(dxdt)}\\
&\leq \,c\,\| e^{\alpha (\frac xR+\varphi(t))^2}(\partial_t+\p_x^3+a_1(x,t)\p_x+a_0(x,t))g\|_{L^2(dxdt)}
\end{aligned}
\end{equation}
holds, for $R\geq R_0$, $\alpha$ such that $\alpha^2\geq M_1R^3$ and $g\in C_0^\infty(\rr^{2})$ supported  in 
$$
\{(x,t)\in\rr^2 : |\tfrac xR+\varphi(t)|\ge 1\}.
$$

\end{lemma}
\begin{proof}
\label{Proof of Lemma 3.2}
From \eqref{3.1}, Lemma 3.1 it follows that
\begin{equation}
\label{3.13.b}
\begin{aligned}
&\frac{\alpha^{5/2}}{R^3}\, \|e^{\alpha (\frac xR+\varphi(t))^2} (\tfrac xR+\varphi(t))^2g\|_{L^2(dxdt)}
+\frac{\alpha^{3/2}}{R^2}\, \|e^{\alpha (\frac xR+\varphi(t))^2} \p_xg\|_{L^2(dxdt)}\\
&\leq c\,\| e^{\alpha (\frac xR+\varphi(t))^2}(\partial_t+\p_x^3)g\|_{L^2(dxdt)}\\
&\leq c \,\| e^{\alpha (\frac xR+\varphi(t))^2}(\partial_t+\p_x^3+a_1(x,t)\p_x+a_0(x,t))g\|_{L^2(dxdt)}\\
&+c \,\| e^{\alpha (\frac xR+\varphi(t))^2}(a_1(x,t)\p_x+a_0(x,t))g\|_{L^2(dxdt)}\\
&\leq c \,\| e^{\alpha (\frac xR+\varphi(t))^2}(\partial_t+\p_x^3+a_1(x,t)\p_x+a_0(x,t))g\|_{L^2(dxdt)}\\
&+ c \|a_1\|_{L^{\infty}_{xt}}\,\| e^{\alpha (\frac xR+\varphi(t))^2} \p_x g\|_{L^2(dxdt)}
+c \|a_0\|_{L^{\infty}_{xt}}\,\| e^{\alpha (\frac xR+\varphi(t))^2}  g\|_{L^2(dxdt)}.
\end{aligned}
\end{equation}

 Since our hypothesis guarantee  that 
$\alpha^{5/2}/R^3$ and $\alpha^{3/2}/R^2$  growth as a positive (fractional) power of  $R$
for $R$ sufficiently large the last two terms in the right hand side of \eqref{3.13.b}  can be hidden in the left hand side to obtain the desired result.

\end{proof}

 \begin{theorem}
\label{Theorem 3.3}
Let $u\in C([0,1]:H^2(\rr))$ be a solution of
\begin{equation}
\label{3.14}
\p_t u+\p_x^3 u +a_2(x,t)\p_x^2 u + a_1(x,t) \p_x u +a_0(x,t) u =0,
\end{equation} 
with $a_0,\,a_1, \,a_2,\, \p_xa_2, \,\p_x^2 a_2\in L^{\infty}(\mathbb R^2)$ 
and  $a_2,\, \p_ta_2\in L_t^{\infty}(\mathbb R:L^1_x(\mathbb R))$. Assume that
$$
\int_{\mathbb R}\int_0^1 (u^2+(\p_x u)^2 + (\p_x^2 u)^2)(x,t)dx dt\leq A^2,
$$
and 
$$
\int_{1/2-1/8}^{1/2+1/8}\,\int_{0<x<1} u^2(x,t)dx dt \geq 1.
$$
Then there exist constants $R_0,\,c_0,\,c_1>0$ depending on
\begin{equation}
\label{constants}
A,\,\|a_0\|_{\infty},\,\|a_1\|_{\infty},\,\|a_2\|_{\infty},\,\|\p_x a_2\|_{\infty},\,\|\p_x^2 a_2\|_{\infty},\,
\|a_2\|_{L^{\infty}_tL^1_x},\,\text{and}\; \;\|\p_t a_2\|_{L^{\infty}_tL^1_x}
\end{equation}
 such that for $R\geq R_0$ 
\begin{equation}
\label{3.15}
\delta(R)=\delta_u(R)=\left(\int_0^1\int_{R-1<x<R}(u^2+(\p_x u)^2 + (\p_x^2 u)^2)dx dt\right)^{1/2}\geq c_0
e^{-c_1R^{3/2}}.
\end{equation}

\end{theorem}

\begin{proof}
\label{Proof of Theorem 3.3} First,  we use a gauge transformation (i.e.  a change of the dependent variable)  to reduce the equation in \eqref{3.14}
 to an  \lq\lq equivalent" one  which does not  involve second order derivative. Define
\begin{equation}
\label{3.15b}
v(x,t)=u(x,t) e^{\frac13\int_0^x\,a_2(s,t)ds}.
\end{equation}
Thus multiplying the equation in \eqref{3.14} by $ e^{\frac13\int_0^x\,a_2(s,t)ds}$ and using that
$$
\aligned
 &e^{\frac13\int_0^x\,a_2(s,t)ds}\,\p_t u = \partial_t v - \tfrac13(\int_0^x\p_ta_2(s,t)ds)v,\\
& e^{\frac13\int_0^x\,a_2(s,t)ds}\,\p_x u = \p_xv-\tfrac13 a_2v,\\
& e^{\frac13\int_0^x\,a_2(s,t)ds}\,\p_x^2u =\p_x^2v-\tfrac23 a_2 \p_xv +
(-\tfrac13\p_xa_2+(\tfrac13a_2)^2) v,\\
 &e^{\frac13\int_0^x\,a_2(s,t)ds}\,\p_x^3u=\p_x^3v-a_2\p_x^2v\\
&\;\;\;\;+((\tfrac13a_2)^2-\tfrac13\p_xa_2)\p_xv +(-(\tfrac13 a_2)^3+\tfrac13a_2\p_xa_2-\tfrac13\p_x^2a_2)v.
\endaligned
$$
the equation for $v=v(x,t)$ can be written as
$$
\p_t v +\p_x^3 v + \tilde a_1(x,t) \p_x v +\tilde a_0(x,t) v =0,
$$
where from our hypothesis on $\,a_0,\,a_1,\,a_2\,$ it follows that $\tilde a_1,\,\tilde a_0\in L^{\infty}(\mathbb R^2)$.

Next, we shall follow the arguments in \cite{EsKePoVe}.
 
For $R>2$ let $\theta_R\in C^{\infty}(\rr)$ with $\theta_R(x)=1$ if $x<R-1$, $\theta_R(x)=0$ if $x>R$.

 Let $\mu\in C^{\infty}(\rr)$ with $\mu(x)=0$ if $x<1$ and $\mu(x)=1$ if $x>2$, and $\varphi\in C^{\infty}_0(\rr)$, $\varphi:\mathbb R\rightarrow [0,3]$ with  
 \begin{equation}
 \label{1234}
 \varphi(t) = 
 \begin{cases}
 \begin{aligned}
 &0,\;\;\;\;t\in[0,1/4]\cap[3/4,1],\\
 &3,\;\;\;\; t\in[1/2-1/8,1/2+1/8].
 \end{aligned}
 \end{cases}
 \end{equation}
 
 We define  
 \begin{equation}
 \label{3.15c}
 g(x,t)=\theta_R(x)\,\mu(\tfrac xR+\varphi(t))\,v(x,t),\;\;\;\;\;\;\;(x,t) \in\mathbb R\times [0,1]
 \end{equation}
and observe that
 
$\cdot$ if $x>R$, then $g(x,t)=0$,

$\cdot$ if $x<R$ and $t\in[0,1/4]\cap[3/4,1]$, then  $g(x,t)=0$,

$\cdot$ if $\tfrac xR+\varphi(t)<1$,  then $g(x,t)=0$, so that $g$ has support on $\mathbb R\times (0,1)$ and can be assumed to satisfy the hypothesis of Lemma 3.1. 

Also, for $(x,t)\in(0,R-1)\times [1/2-1/8,1/2+1/8]$, $g(x,t)=v(x,t)$ and $|\tfrac xR+\varphi(t)|\geq 2$. 
 
 From \eqref{3.15c} one has that
 \begin{equation}
\label{3.16}
\begin{aligned}
&(\partial_t +\p_x^3+\tilde a_1\p_x+\tilde a_0)g\\
&=\mu(\tfrac xR+\varphi(t))\left[3\theta^{(1)}_R\p_x^2v +3\theta^{(2)}_R\p_xv +\theta^{(3)}_Rv +\tilde a_1\theta^{(1)}_Rv\right]\\
&+\theta_R(x)\left(\mu^{(1)}(\cdot)(\varphi^{(1)}+\tfrac{\tilde a_1}{R})v
+3\mu^{(1)}(\cdot)\tfrac1R\p_x^2v
+3\mu^{(2)}(\cdot)\tfrac{1}{R^2}\p_xv
+\mu^{(3)}(\cdot)\tfrac{1}{R^3}v\right)\\
&+3\mu^{(1)}(\cdot)\tfrac{1}{R}\theta^{(2)}_R v +3\mu^{(2)}(\cdot)\tfrac{1}{R^2}\theta^{(1)}_R v+6\mu^{(1)}(\cdot)\tfrac{1}{R}\theta^{(1)}_R\partial_x v,
\end{aligned}
\end{equation}
 where the first term in the right hand side of \eqref{3.16} is supported in
 $[R-1,R]\times[0,1]$, where $|\tfrac xR+\varphi(t)|\leq 4$, and the remaining terms in the right hand side 
 of \eqref{3.16} are supported in $\{(x,t)\,:\,1\leq |\tfrac xR+\varphi(t)|\leq 2\}$.
 
 Using the notation
  \begin{equation}
\label{3.16.b}
 \delta_v(R)=\left(\int_0^1\int_{R-1<x<R}(v^2+(\p_x v)^2 + (\p_x^2 v)^2)(x,t)dx dt\right)^{1/2},
\end{equation}
 from \eqref{3.13} and \eqref{3.16} it follows that
 $$
 c\frac{\alpha^{5/2}}{R^3}e^{4\alpha}\leq c\frac{\alpha^{5/2}}{R^3}\|
 e^{\alpha(\tfrac xR+\varphi(t))^2}g\|_{L^2(dxdt)} \leq c_1e^{16\alpha}\,\delta_v(R)+c_2e^{4\alpha}A,
$$
therefore
$$
 c\frac{\alpha^{5/2}}{R^3}\leq c_1e^{12\alpha}\,\delta_v(R) +c_2A.
 $$
 
 Taking $\alpha=M_1R^{3/2}$ with $M_1$ as in Lemma 3.2 we get
 \begin{equation}
 \label{3.17}
 cM_1^{5/2} R^{3/4}\leq c_1e^{12M_1R^{3/2}}\,\delta_v(R)+c_2 A.
 \end{equation}

For $R$ sufficiently large the last term in the right hand side of \eqref{3.17} can be absorbed 
into the 
left hand side to get that
$$
\delta_v(R) \geq \frac c2 M_1^{5/2} R^{3/4}e^{-12M_1R^{3/2}}.
$$
 Finally, from  \eqref{3.15b}, \eqref{3.15c}, \eqref{3.16.b}  and our hypothesis one has that $\delta_u\sim \delta_v$, i.e. there exists
 $c>1$ such that 
 $$
 c^{-1}\delta_v(R)\leq \delta_u(R)\leq c \delta_v(R),\;\;\;\;\;\forall\, R\geq R_0,
 $$
 which yields the desired result.
 \end{proof}

 \section{Proof of Theorems 1.3 and 1.1}
   
  \it{Proof of Theorem 1.3. }\rm

  If $u \not \equiv 0$, we can assume, after a possible translation, dilation,  and multiplication by a constant, that $u=u(x,t)$ satisfies the  hypothesis of Theorem 3.1. Hence, we have that \begin{equation}
\label{4.1}
 \delta_u(R)=\left(\int_0^1\int_{R-1<x<R}(u^2+(\p_x u)^2 + (\p_x^2 u)^2)(x,t)dx dt\right)^{1/2}
 \geq c_0\,e^{-c_1R^{3/2}},
\end{equation}
for all $R$ sufficiently large where the constants $c_0, \,c_1$ depend on the quantities in \eqref{constants}.

Now we apply Theorem 2.1 with $\alpha=3/2$ and $a>>4
^{3/2}c_1$ with $c_1$ as in \eqref{4.1}
to conclude that
\begin{equation}
\label{4.2}
 \delta_u(R)\leq c e^{-aR^{3/2}/4^{3/2}},
 \end{equation}
for all $R$ sufficiently large. Combining \eqref{4.1} and \eqref{4.2} and letting $R \uparrow \infty$ we get a contradiction.
Therefore $u\equiv 0$.
 

\vskip.1in

\it{Proof of Theorem 1.1. }\rm 

It will be shown that Theorem 1.3 applies to the equation of the type \eqref{1.2} satisfied by the difference
$u_1-u_2$ of the solutions. Thus, one just needs to prove that the coefficients
$a_0,\,a_1$ satisfy the assumptions \eqref{1.10} and \eqref{1.10.c} . We recall that in this case $a_2\equiv 0$.

 Since for any $k\in \zz^+$,  $a_0,\,a_1$ are  polynomials of order $k$ in $u_1,\,u_2, \,\p_xu_1$,  
 with $u_1,\,u_2\in C([0,1]:H^2(\rr))$, and $a_2\equiv 0$ it is clear that the hypothesis \eqref{1.10.c}
 holds. So it remains to check the 
 conditions \eqref{1.10}, i. e.
 \begin{equation}
 \label{5.5.5}
\begin{aligned}
&a_0\in L^{4/3}_{xt}\cap L_x^{16/13}L^{16/9}_t\cap L_x^{8/7}L_t^{8/3},\\
&a_1\in L^{16/13}_xL^{16/9}_t\cap L_x^{8/7}L^{8/3}_t\cap L_x^{16/15}L_t^{16/3}.\\
\end{aligned}
\end{equation}

First, we consider the KdV equation, i. e. $k=1$ in \eqref{1.1}, for which we have
$$
a_0(x,t)=\p_x u_1(x,t)\;\;\;\;\;\;\text{and}\;\;\;\;\;\;a_1(x,t)= u_2(x,t).
$$

Using the hypothesis   $u_1, \,u_2\in C([0,1] : H^3 \cap L^2(|x|^2dx)$
it follows by interpolation (or integration by parts) that
\begin{equation}
\label{5.1}
a_0\in L^{\infty}([0,1] : H^2)\;\;\;\;\;\;\text{and}\;\;\;\;\;\;|x|^{2/3}a_0,\;|x|^{1/3}\p_x a_0\in L^{\infty}([0,1] : L^2_x),
\end{equation}
and by Sobolev lemma that
\begin{equation}
\label{5.2}
|x|^{1/3} \, a_0\in L^{\infty}([0,1] : L_x^{\infty}).
\end{equation}

Thus,  \eqref{5.1} and H\"older inequality yields 
$$
\|a_0(t)\|_{L^{4/3}_{xt}}\leq c \sup_{t\in[0,1]}\, \|(1+|x|)^{1/2}\,a_0(\cdot, t)\|_{L^2},\;\;\;\;\;t\in\rr,
$$
which proves that $a_0\in L^{4/3}_{xt}$.  Next, the string of inequalities,
\begin{equation}
\label{5.3}
\begin{aligned}
&\|a_0\|_{L_x^{16/13}L^{16/9}_t}\\
&=\left(\int \frac{1}{(1+|x|)^{8/13}}\,(1+|x|)^{8/13}\,
\left(\int |a_0(x,t)|^{16/9}dt\right)^{9/13}dx\right)^{13/16}\\
&\leq c\left(\int\,\| (1+|x|)^{1/2} \, a_0(\cdot,t)\|_{L^2_x}^{16/9}\,dt\right)^{9/16}\\
&\leq \sup_{t\in[0,1]}\,(\| a_0(\cdot,t)\|_{L^2_x}+\| |x|^{1/2} \,a_0(\cdot,t)\|_{L^2_x}),
\end{aligned}
\end{equation}
and \eqref{5.1}  show that $\,a_0\in L_x^{16/13}L^{16/9}_t$.

In a similar fashion we have that 
\begin{equation}
\label{5.4}
\begin{aligned}
&\|a_0\|_{L_x^{8/7}L^{8/3}_t}\\
&=\left(\,\int \,\frac{1}{(1+|x|)^{4^+/7}}\,(1+|x|)^{4^+/7}\,
\left(\int |a_0(x,t)|^{8/3}dt\right)^{3/7}dx\right)^{7/8}\\
&\leq c \left(\int \int (1+|x|)^{4^+/3}|a_0(x,t)|^{8/3}dx dt\right)^{3/8}\\
&\leq c \|(1+|x|^{2/3}) \,a_0\|_{L^2_{xt}}^{1/3}\, \|(1+|x|^{\epsilon}) \,a_0\|_{L^{\infty}_{tx}}^{2/3},
\end{aligned}
\end{equation}
for any $\epsilon>0$ which together with  \eqref{5.1} and \eqref{5.2}  imply that $\,a_0\in L_x^{8/7}L^{8/3}_t$.

Now we consider $a_1(x,t)=u_2\in C([0,1] : H^3 \cap L^2(|x|^2dx)$. Thus, it follows that
\begin{equation}
\label{5.5}
|x| a_1,\;|x|^{2/3}\p_x a_1,\;|x|^{1/3}\p^2_x a_1,\;\p_x^3 a_1 \in L^{\infty}([0,1] : L^2_x),
\end{equation}
and  by Sobolev lemma  that
\begin{equation}
\label{5.6}
|x|^{2/3} a_1\in L^{\infty}([0,1] : L_x^{\infty}).
\end{equation}

The same arguments used in \eqref{5.3} and \eqref{5.4} show  that 
$a_1 \in L_x^{16/13}L^{16/9}_t\cap L_x^{8/7}L^{8/3}_t$. So it only remains to prove  that
$a_1\in L_x^{16/15}L_t^{16/3}$. A familiar process  leads to
\begin{equation}
\label{5.7}
\begin{aligned}
&\|a_0\|_{L_x^{16/5}L^{16/3}_t}\\
&\leq \left(\int  \int (1+|x|)^{4^+} |a_1(x,t)|^{16/3}dx dt\right)^{3/16}\\
&\leq c \left(\int \int (1+|x|^2)\,|a_1(x,t)|^{2}dx dt\right)^{3/16}\
\, \|(1+|x|^{2^+}) |a_1|^{10/3}\|_{L^{\infty}_{xt}}^{3/16}\\
&\leq c\,\|(1+|x|) a_1\|^{3/8}_{L^2_{xt}}\,\||x|^{3^+/5}\,a_1\|^{5/8}_{L^{\infty}_{xt}} \leq c\,\|(1+|x|) a_1\|^{3/8}_{L^2_{xt}}\,\||x|^{2/3}\,a_1\|^{5/8}_{L^{\infty}_{xt}}.
\end{aligned}
\end{equation}
Therefore, inserting \eqref{5.5}, \eqref{5.6}  in \eqref{5.7} one obtains the desired result. 

We have completed the proof of Theorem 1.1. in the case of the KdV equation.

Next, we turn to the proof of Theorem 1.1 for the equations in \eqref{1.1} with $k\geq 2$.
Using that $u_1,\,u_2\in L^{\infty}(\rr\times [0,1])$ it suffices to consider the case $k=2$ where 
$$
a_0(x,t)=(u_1+u_2)\p_xu_1,\;\;\;\;\;\;\;\;\;a_1(x,t)=u_2^2.
$$ 
Since  $\,u_1,\,u_2\in C([0,1] : H^2\cap L^2(|x|^2dx))$ by interpolation and Sobolev lemma it follows that
$$
|x| u_j,\;\;\;\;|x|^{1/2}\p_xu_j \in L^{\infty}([0,1] : L^2_x),\;\;\;\;|x|^{1/2}u_j\in L^{\infty}([0,1] : L_x^{\infty}),\;\;\;\;j=1,2.
$$
Hence,
$$
|x|^{3/2}a_0\in L^{\infty}([0,1] : L^2_x),\;\;|x|^{1/2}a_0\in  L^{\infty}([0,1] : L_x^{\infty}),
$$
and 
$$
|x|a_1\in L^{\infty}([0,1] : L^2_x),\;\;|x|^{2/3}a_1\in  L^{\infty}([0,1] : L_x^{\infty}),
$$
which were the conditions used  to obtain the result in the case $k=1$.

This completes the proof of Theorem 1.1.

\vspace{3mm} \noindent{\large {\bf Acknowledgments}}
\vspace{3mm}\\  L. E. and L. V. were supported by a MEC grant and by the European Comission via the network Harmonic Analysis and Related Problems. C. E. K. and G. P. were supported by NSF grants.


\vskip1cm

\noindent{\bf Luis Escauriaza}\\
Departamento de Matematicas\\
Universidad del Pais Vasco\\
Apartado 644\\
48080 Bilbao\\
Spain\\
E-mail: mtpeszul@lg.ehu.es\\

\noindent{\bf Carlos E. Kenig}\\
Department of Mathematics\\
University of Chicago\\
Chicago, Il. 60637 \\
USA\\
E-mail: cek@math.uchicago.edu\\

\noindent{\bf Gustavo Ponce}\\
Department of Mathematics\\
University of California\\
Santa Barbara, CA 93106\\
USA\\
E-mail: ponce@math.ucsb.edu\\

\noindent{\bf Luis Vega}\\
Departamento de Matematicas\\
Universidad del Pais Vasco\\
Apartado 644\\
48080 Bilbao\\
Spain\\
E-mail: mtpvegol@lg.ehu.es

\end{document}